\documentclass{amsart}

\usepackage{amsmath,amsthm,amssymb,amscd,mathabx,accents,enumitem,mathtools}
\usepackage{amsfonts}
\usepackage{rotating}
\usepackage{euscript}
\usepackage{pst-node}
\usepackage{epsfig,verbatim}
\usepackage{tikz-cd}
\def\be{\begin{equation}}
\def\ee{\end{equation}}

\def\C{{\mathbb C}} 
\def\f{\EuScript}
\def\N{{\mathbb N}} 
\def\P{{\mathbb P}}

\def\phi{{\varphi}}
\def\v{{\varepsilon}} 

\def\deg{{\rm deg\,}}

\def\Aut{{\rm Aut\,}}

\def\Ord{{\rm Ord}}

\def\bp{\begin{proposition}}
\def\ep{\end{proposition}}

\def\bt{\begin{theorem}}
\def\et{\end{theorem}}
\def\br{\begin{remark}}
\def\er{\end{remark}}
\def\be{\begin{equation}}
\def\bee{\begin{equation*}}
\def\l{\label}

\def\ee{\end{equation}}
\def\eee{\end{equation*}}
\def\bl{\begin{lemma}}
\def\el{\end{lemma}}
\def\bc{\begin{corollary}}
\def\ec{\end{corollary}}
\def\pr{\noindent{\it Proof. }}

\def\bd{\begin{definition}}
\def\ed{\end{definition}}

\mathtoolsset{showonlyrefs}

\newtheorem{theorem}{Theorem}[section]
\newtheorem{lemma}[theorem]{Lemma}
\newtheorem{corollary}[theorem]{Corollary}
\newtheorem{proposition}[theorem]{Proposition}
\theoremstyle{definition}
\newtheorem{remark}[theorem]{Remark}

\begin{document}
\title{Sharing a measure of maximal entropy 
in polynomial semigroups
}
\author{Fedor Pakovich}
\thanks{
This research was supported by ISF Grant No. 1432/18}
\address{Department of Mathematics, Ben Gurion University of the Negev, Israel}
\email{
pakovich@math.bgu.ac.il}

\begin{abstract}
Let $P_1,P_2,\dots, P_k$ be complex polynomials of degree at least two that are 
not simultaneously conjugate to monomials or to Chebyshev polynomials, and $S$ the semigroup  under composition  generated by $P_1,P_2,\dots, P_k$.  We show that 
all elements of $S$ share a measure of maximal entropy 
if and only if the intersection of principal right ideals  $SP_1\cap SP_2\cap \dots \cap SP_k$ is non-empty.
\end{abstract}

\maketitle

\section{Introduction}
In the recent paper by Jiang and Zieve \cite{zi},  the authors showed that a semigroup  of polynomials    under composition generated by two complex polynomials $P_1$ and $ P_2$ of degrees $n_1\geq 2$ and $n_2\geq 2$   
is not free if and only if it is isomorphic to the semigroup generated by $z^{n_1}$ and $\v z^{n_2}$, where 
 $\v$ is a root of unity. This implies in particular that whenever $S=\langle P_1,P_2\rangle$ is not free  
there exists $r>0$ for which  $P_1^{\circ r}$ and $P_1^{\circ r}\circ P_2$
commute. 
Since commuting polynomials can be described explicitly (\cite{r3},\cite{r2}), the last property permits to classify all pairs of polynomials $P_1$ and $P_2$ for which $S=\langle P_1,P_2\rangle$ is not free. 

Being combined with the description of pairs of rational functions sharing a measure of maximal entropy obtained in \cite{lev}, \cite{lp},  the result of \cite{zi}   implies  the following criterion: a semigroup $S=\langle P_1,P_2\rangle$ generated by two  polynomials $P_1$ and $P_2$ of degree at least two that are not simultaneously conjugate to monomials or to Chebyshev polynomials is not free if and only if all elements of $S$ 
share a measure of maximal entropy. The problem of characterization of  semigroups of polynomials satisfying the last property has been studied in the recent papers \cite{peter}, \cite{p20}, where several equivalent characterizations of such semigroups in semigroups-theoretic terms were given. Among such characterizations we mention the right amenability and the absence of free sub\-semi\-groups. The result of \cite{zi} provides yet another  characterization of such semigroups,  in terms of freeness, in case if  considered semigroups  are generated by  two  polynomials.

It is not hard to see that  the result of \cite{zi} and the corresponding characterization  of semigroups whose elements share a measure of maximal entropy do not allow for a direct generalization to a greater number of generators. For example, for arbitrary polynomials $R$ and $X$ setting 
$$P_1=R\circ z^n, \ \ \ P_2=X, \ \ \ P_3=\v X,$$ where $\v$ satisfies $\v^n=1$, we 
obtain a semigroup $S=\langle P_1,P_2,P_3\rangle$, which is not free  since 
 $$P_1\circ P_2=P_1\circ P_3.$$ However, it is clear that in general $P_1,P_2,P_3$ do not  share a measure of maximal entropy.

In this note, we provide a generalization of the result of \cite{zi} 
to arbitrary finitely generated semigroups of polynomials  replacing the non-freeness condition by an other condition, which is however equivalent to the condition that $S$ is not free in case if the number of generators equals two. We also provide a  
characterization  of finitely generated semigroups of polynomials  whose elements share a measure of maximal entropy. 

 To formulate our results explicitly, we introduce some notation. 
We say that two semigroups of polynomials $S_1$ and $S_2$ are conjugate if 
there exists $\alpha\in \Aut(\C)$ such that $$\alpha\circ S_1\circ \alpha^{-1}=S_2.$$ We  denote by $\f Z$ the  semigroup of polynomials  consisting of  monomials  $az^n,$ where $a\in \C^*$ and $n\geq 1,$ and by $\f T$ the   semigroup consisting of  polynomials of the form $\pm T_n,$ $n\geq 1,$ where $T_n$ stands for the  Chebyshev polynomial of degree $n$.  
Finally, we denote  by  $\f Z^U$ the subsemigroup of $\f Z$ consisting of  polynomials  of 
the form $\omega z^n,$ where $\omega$ is a root of unity.

In this notation, our first result is following.

\bt \l{main0} Let $P_1,P_2,\dots, P_k$ be complex polynomials of degree at least two.
Then the semigroup $S=\langle P_1,P_2,\dots, P_k\rangle$ is isomorphic to a subsemigroup of $\f Z^U$ 
if and only if the intersection of principal right ideals    $SP_1\cap SP_2\cap \dots \cap SP_k$ is non-empty.
\et

It is easy to see that for $k=2$ the condition 
\be \l{int} SP_1\cap SP_2\cap \dots \cap SP_k\neq \emptyset\ee
is equivalent to the condition that $S$ is not free. Indeed, any semigroup of rational functions is right cancellative. Therefore, if there exist two different words in the letters $\{P_1,P_2\}$ representing the same element in $S=\langle P_1,P_2\rangle,$ then cancelling  their common suffix we obtain 
two different words representing the same element with different ending letters. Since one of these letters is $P_1$ and the other one is $P_2$, this implies  that $SP_1\cap SP_2\neq \emptyset.$
However, for $k>2$ condition \eqref{int}  
is clearly stronger than merely the requirement that  $S$ is not free.

Our second result is following.

\bt \l{main} Let $P_1,P_2,\dots, P_k$ be complex polynomials of degree at least two such that  $S=\langle P_1,P_2,\dots, P_k \rangle$ is not conjugate to a subsemigroup of $\f Z$ or $\f T$. Then  all elements of $S$ share a measure of maximal entropy if and only if the intersection of principal right ideals  $SP_1\cap SP_2\cap \dots \cap SP_k$  is non-empty.
\et

Notice that since in the polynomial case having the same measure of maximal entropy is equivalent to having the same Julia set,  
Theorem \ref{main} can be viewed as a characterization of polynomials $P_1,P_2,\dots, P_n$ sharing a  Julia set 
via existence of a relation of the form 
\be \l{rel} A_1 P_1= A_2 P_2=\dots = A_n P_n,\ee
where $A_i,$ $1\leq i \leq n,$ are words in $P_1,P_2,\dots P_n$.

The assumption that $S$  is not conjugate to a subsemigroup of $\f Z$ or $\f T$ is not essential for ``if'' part of Theorem \ref{main},  but essential for ``only if'' part. Indeed,  for instance, polynomials $z^n$ and  $bz^m$, $b\in \C^*$, share a measure of maximal entropy whenever $\vert b \vert =1$, but generate a free group, unless $b$ is a root of unity.     

Finally, notice that since semigroups of polynomials whose elements share a mea\-sure of maximal entropy admit many equivalent descriptions (see \cite{p20}), Theorem \ref{main} also can be formulated in many equivalent forms. In particular, under assumptions of Theorem \ref{main}, condition \eqref{int} 
is equivalent to the condition that there exists a 
polynomial $T$ of the form $T=z^rR(z^l),$ where $R\in \C[z]$, $l\geq 1$, and $0\leq r<l$, 
such that 
$$P_i=\omega_iT^{\circ l_i},\ \ \ \ 1\leq i \leq k,$$
for some $l_i$, $1\leq i \leq k,$  and  $l$th roots of unity $\omega_i,$ $1\leq i \leq k.$

\section{Proof of Theorem \ref{main0}} 
Let us recall that for every complex polynomial $P_1$ of degree $n_1\geq 2$ there exists a series $$\beta=\sum_{i=-1}^{\infty}c_iz^{-i}, \ \ \ \ \ \ \  c_{-1}\neq 0,$$ called a B\"otcher function, which makes the diagram  
\be \l{eg} 
\begin{CD} 
\C\P^1 @>P_1>> \C\P^1 \\
@VV \beta V @VV \beta V\\ 
\C\P^1 @>z^{n_1}>> \C\P^1\
\end{CD}
\ee
commutative. Having its roots in complex dynamics (see \cite{bo}, \cite{mil}), B\"otcher function is widely used for studying functional relations between polynomials and related problems (see \cite{a2}, \cite{b}, \cite{be2}, \cite{be1}, \cite{dw}, \cite{zi}, \cite{r3}, \cite{sh}).

As in the paper \cite{zi}, our proof of the ``if'' part of Theo\-rem \ref{main0} uses B\"otcher function and the following lemma  proved in the paper \cite{dw}. 
Following \cite{dw}, for a non-zero element $T=b_0+b_1z+b_2t^2+\dots $ of $ \C[[z]]$,  we define $\Ord_0(T)$ as the minimum number $i\geq 0$ such that $b_i\neq 0$, and  
$l_0(T)$ as the difference $m-\Ord_0(T)$, where $m$ is the minimum number greater than $\Ord_0(T)$ such that $b_m\neq 0$. In case if $T$ is a monomial, we set $l_0(T)=\infty.$ The parameter $l_0(T)$ possesses certain  properties  making it useful for studying functional relations between powers series (see \cite{dw}, Lemma 2.6). Below we need only the properties  listed in the following statement,   which can be checked by a direct calculation.

\bl \l{do} Let $X$ be an element of $z\C [[z]]$ such that $l_0(X)<\infty.$ Then for any 
element $T$ of $z\C [[z]]$ with  $l_0(T)<\infty$ the inequality  
\be \l{bu1} l_0(T\circ X)\geq \min\big(l_0(X),\Ord_0(X)l_0(T)\big)\ee 
holds, and the equality is attained whenever  
$l_0(X)\neq \Ord_0(X)l_0(T).$  On the other hand, if $l_0(T)=\infty$, then the equalities  
\be \l{bu2} l_0(X\circ T)=\Ord_0(T)l_0(X),\ee 
\be \l{bu3}  l_0(T\circ X )=l_0(X)\ee 
hold. \qed
\el

Lemma \ref{do} implies the following corollary. 

\bc \l{kora} Let $T_i,$ $1\leq i \leq k,$  be elements of $z^2\C [[z]]$ such that $l_0(T_1)=l<\infty,$ and  \be \l{ti} l_0(T_{1})\leq l_0(T_{i}), \ \ \ \ 2\leq i \leq k.\ee Assume that  
 $A$ is a word in $T_1,T_2,\dots T_k$, and  $X$ is an element of $z^2\C [[z]]$. Then 
$ l_0(AX)=l,$ 
if $l_0(X)=l$, and 
$ l_0(AX)>l,$ 
if $l_0(X)>l$.
\ec
\pr In the both cases, the proof is by induction on the length of  $A$. If $A$ is empty, then the corollary is trivially true. Further, in the first case, the induction step reduces to the following statement: if $X\in z^2\C [[z]]$ satisfies $l_0(X)=l$, then $$l_0(T_iX)=l, \ \ \ \ 1\leq i \leq k.$$ In turn, the last statement follows from formulas 
\eqref{bu1}, \eqref{bu3} taking into account that inequality \eqref{ti} implies  the inequality 
\be \l{gope} \Ord_0(X)l_0(T_i)\geq \Ord_0(X)l_0(T_1)\geq  2l>l, \ \ \ 1\leq i \leq k.\ee

Similarly, in the second case, we must prove that if $l_0(X)>l$, then \be \l{gav} l_0(T_iX)>l, \ \ \ \  1\leq i \leq k.\ee If $l_0(X)<\infty$, then \eqref{gav} follows from \eqref{bu1} and \eqref{bu3} taking into account the inequalities $l_0(X)>l$ and \eqref{gope}.  
 On the other hand, if $l_0(X)=\infty$, then either $l_0(T_i)=\infty$ and 
 $l_0(T_iX)=\infty >l,$ or $l_0(T_i)<\infty$ and 
$$l_0(T_iX)=\Ord_0(X)l_0(T_i)>l$$
by \eqref{bu2} and \eqref{gope}. \qed

We deduce Theorem \ref{main0} from the following result. 

\bt \l{main1}   Let  $Q_i,$ $1\leq i \leq k,$ be elements of $z^2\C [[z]]$, and $R$ the semigroup  generated by $Q_1, Q_2,\dots, Q_k$. Assume that $Q_1$ is contained in $\f Z^U$. Then  
\be \l{xera} RQ_1\cap RQ_2\cap \dots \cap RQ_k\neq\emptyset\ee
if and only if every $Q_i,$  $2\leq i \leq k$, is contained in $\f Z^U$.
 \et
\pr 
 Let  $Q_i=a_{i,1}z^{n_i}+ a_{i,2}z^{n_i+1}+\  \dots\ , $ $1\leq i \leq k,$ where $a_{i,1}\neq 0.$ 
Assume that \eqref{xera} holds, but not all $Q_i,$ $1\leq i \leq k,$ are monomials. 
Without loss of generality we may assume that for some $s$, $1\leq  s<k,$ the series $Q_1,Q_2,\dots , Q_s$ are monomials, while the series $Q_{s+1},\dots, Q_k$ are not, and that 
\be \l{ta} l_0(Q_{s+1})\leq \dots \leq l_0(Q_k).\ee
By condition, there exist words $A_1,A_2,\dots, A_k$ in  $Q_1,Q_2,\dots , Q_k$ 
such that 
\be \l{rell} A_1Q_1=A_2Q_2=\dots =A_{k}Q_{k}. \ee 
Applying the first part of Corollary \ref{kora} to the word $A_{s+1}Q_{s+1},$ we obtain that $$l_0(A_{s+1}Q_{s+1})=l_0(Q_{s+1}).$$  On the other hand, applying the second part of Corollary \ref{kora} to the word $A_{1}Q_{1},$ we obtain that $$l_0(A_{1}Q_{1})>l_0(Q_{s+1}).$$ Since $A_{1}Q_{1}=A_{s+1}Q_{s+1}$, we obtain a contradiction, which shows that all $Q_1,Q_2,$ $\dots , Q_s$ are monomials. In particular, equality \eqref{rell} reduces to  
 the equality \be \l{rel0} A_1a_{1,1}z^{n_1}=A_2a_{2,1}z^{n_2}=\dots =A_ka_{k,1}z^{n_k},\ee
where $A_i,$ $1\leq i \leq k,$ are words in $a_{1,1}z^{n_1},a_{2,1}z^{n_2},\dots ,a_{k,1}z^{n_k}.$

Clearly, \eqref{rel0} implies an equality of the form 
\be \l{krot} U_1z^N=U_2z^N=\dots =U_kz^N,\ee where $U_i,$ $1\leq i \leq k,$ are monomials 
in $a_{1,1},a_{2,1},\dots , a_{k,1}$, and $N$ is a natural number. To finish the proof of the ``if'' part of 
the theorem it is enough to show that whenever $i\neq j,$ $1\leq i,j \leq k,$ the inequality \be \l{ft} \deg_{a_{i,1}} U_i>\deg_{a_{i,1}} U_j\ee holds. Indeed, in this case making in the equality 
$$U_1=U_2=\dots =U_k$$ all possible cancellations, we obtain an equality of the form 
\be \l{xey} a_{1,1}^{s_1}=a_{2,1}^{s_2}=\dots =a_{k,1}^{s_k},\ee
where $s_i\geq 1,$ $1\leq i \leq k,$ implying that all  $a_{i,1},$  $2\leq i \leq k$,
are roots of unity.

It is clear that 
the minimum  value of $\deg_{a_{i,1}} U_i$ is attained if the word $A_i$ contains no letter $a_{i,1}z^{n_i}$ at all, implying that  
$$\deg_{a_{i,1}} U_i\geq \frac{N}{n_i}.$$ 
Thus, to prove \eqref{ft} it is enough to show that 
\be \l{doi} \deg_{a_{i,1}} U_j< \frac{N}{n_i}.\ee
Let $r$ be the number of appearances of $a_{i,1}z^{n_i}$ in 
 $A_j$. It is easy to see that the maximum value of $\deg_{a_{i,1}} U_j$ 
is attained if these appearances occur in the last $r$ letters of   
 $A_j$, implying that  
\be \l{df} \deg_{a_{i,1}} U_j\leq \frac{N}{n_jn_i}+ \frac{N}{n_jn_i^2}+\frac{N}{n_jn_i^3}+\dots +\frac{N}{n_jn_i^r}.\ee  Taking into account that $n_i\geq 2,$ $n_j\geq 2$, this implies that 
$$\deg_{a_{i,1}} U_j<\frac{N}{n_jn_i}\sum_{l=0}^{\infty}\frac{1}{n_i^l}\leq \frac{N}{2n_i}\frac{1}{1-\frac{1}{n_i}}\leq \frac{N}{n_i}.$$

Let us assume now that  $Q_i,$  $1\leq i \leq k$, are contained in $\f Z^U,$ and show that then \eqref{xera} holds. 
Let $l\geq 1$ be a number such that all $a_{i,1},$ $1\leq i \leq k$, are $l$th roots of unity. Setting 
$$
F_1=Q_1\circ Q_2\circ \dots \circ Q_k, \ \ \ F_2=Q_2\circ Q_3\circ \dots \circ Q_1, \ \  \dots \ \ , F_k=Q_k\circ Q_1\circ \dots \circ Q_{k-1} 
$$
and observing that $$\deg F_1=\deg F_2=\dots =\deg F_k,$$ we see that for every $j\geq 1$ there exists an $l$th root of unity $\omega_j$ such that 
$$F_1^{\circ j}=\omega_jF_2^{\circ j}.$$  
The pigeonhole principle yields that there exist an infinite subset $K_1$ of $\N$ and an $l$th root of unity $\delta_1$   such that for every $j\in K_1$ the equality  
$$F_1^{\circ j}=\delta_1F_2^{\circ j}$$
holds,  implying that for every $j_1,j_2\in K_1$ with $j_2>j_1$ the equality 
$$F_1^{\circ j_2}=F_1^{\circ j_1}\circ F_2^{\circ (j_2-j_1)}$$ holds. 
Similarly, there exist an infinite subset $K_2$ of $K_1$ and an $l$th root of unity $\delta_2$ such that for every $j\in K_2$ the equality   
$$F_1^{\circ j}=\delta_2F_3^{\circ j}$$
holds, and  for every $j_1,j_2\in K_2$ with $j_2>j_1$ the equality 
$$F_1^{\circ j_2}=F_1^{\circ j_1}\circ F_3^{\circ (j_2-j_1)}$$ holds. 
Continuing in the same way, we will find natural numbers $j_2$ and $j_1$ such that $j_2>j_1$ and 
\be \l{bes} F_1^{\circ j_2}=F_1^{\circ j_1}\circ F_2^{\circ (j_2-j_1)}=F_1^{\circ j_1}\circ F_3^{\circ (j_2-j_1)}=\dots =F_1^{\circ j_1}\circ F_k^{\circ (j_2-j_1)}.\ee Thus,   
$$F_1^{\circ j_2}\in RQ_1\cap RQ_2\cap \dots \cap RQ_k,$$
implying \eqref{xera}. \qed

\vskip 0.2cm
\noindent{\it Proof of Theorem \ref{main0}.}
Let $P_1,P_2,\dots, P_k$ be polynomials of degree at least two. Since the B\"otcher function $\beta$ for $P_1$ 
provides an isomorphism between the semigroup $S=\langle P_1,P_2,\dots, P_k \rangle$ 
and the semigroup of power series $R$  generated by the power series $z^{n_1}$ and 
$$Q_i=\beta \circ P_i\circ \beta^{-1}, \ \ \ \ 2\leq i \leq n,$$ if \eqref{int} holds, then \eqref{xera} also holds implying that 
$S$ is isomorphic to a subsemigroup of $\f Z^U$ by Theorem \ref{main}. 

In the other direction, if $S$ is isomorphic to a subsemigroup of $\f Z^U$, then for the images  $Q_1, Q_2,\dots, Q_k$ of $P_1,P_2,\dots, P_k$ under this isomorphism condition \eqref{xera} holds by Theorem \ref{main}. Therefore, for $P_1,P_2,\dots, P_k$ condition \eqref{int} holds. \qed

\subsection{Proof of Theorem \ref{main}}

Let us recall that if $f$ is a rational function of degree $n\geq 2$, then by the results of Freire, Lopes, Ma\~n\'e (\cite{flm}) and Lyubich (\cite{l}), there exists a unique probability measure $\mu_f$ on $\C\P^1$, which is invariant under $f$, has support equal to the Julia set $J(f)$, and achieves maximal entropy $\log n$ among all $f$-invariant probability measures. It is clear that the equality $\mu_f=\mu_g$ implies the equality of the Julia sets $J(f)=J(g).$ Moreover, for polynomials these conditions are equivalent. The problem of describing rational functions sharing a measure of maximal entropy and the problem of describing rational functions sharing a Julia set have been studied in
\cite{a2}, \cite{b},  \cite{be2}, \cite{be1}, \cite{lev}, \cite{lp}, \cite{p1}, \cite{parn},  \cite{sh}, \cite{ye}. 

For any rational functions $f$ and $g$   the equality
\be \l{zx} f^{\circ j_1}=f^{\circ j_2}\circ g^{\circ s}\ee  for some $j_1,s\geq 1$ and $j_2\geq 0$ 
implies that $f$ and $g$ share the measure of maximal entropy. 
Furthermore, the results of the papers \cite{lev} and \cite{lp} imply that if  the functions $f$ and $g$ are 
neither Latt\`es maps nor conjugate to $z^{\pm n}$ or $\pm T_n,$   
then the equality $\mu_f=\mu_g$ holds if and only if equality \eqref{zx} holds  
(see \cite{parn}, \cite{ye} for more detail).

\vskip 0.2cm
\noindent{\it Proof of Theorem \ref{main}.}
In view of the isomorphism between semigroups $S$ and $R$,  
the proof of the ``if'' part of Theorem \ref{main0} shows that if \eqref{int} holds, then the polynomials 
$$G_1=P_1\circ P_2\circ \dots \circ P_k, \ \ \ G_2=P_2\circ P_3\circ \dots \circ P_1, \ \  \dots \ \ , G_k=P_k\circ P_1\circ \dots \circ P_{k-1}$$ along 
with polynomials $F_i,$ $1\leq i \leq k,$  satisfy relations \eqref{bes}, implying that    
\be \l{xir} J(G_1)=J(G_2)=\dots =J(G_k).\ee 
On the other hand, since the semiconjugacy relation  
\be 
\begin{CD}
\C\P^1@>B>> \C\P^1\\
@VV X V @VV X V\\ 
\C\P^1 @>A>> \C\P^1
\end{CD}
\ee
between rational functions  of degree at least two  implies that 
$$X^{-1}(J(A))=J(B)$$ (see e.g. \cite{buf}, Lemma 5),
the semiconjugacies 
\be 
\begin{CD}
\C\P^1@>G_{1}>> \C\P^1\\
@VV P_k V @VV  P_k V\\ 
\C\P^1 @>G_{k}>> \C\P^1\,,
\end{CD}
\ \ \ \ \ 
\begin{CD}
\C\P^1@>G_{i+1}>> \C\P^1\\
@VV P_i V @VV  P_i V\\ 
\C\P^1 @>G_{i}>> \C\P^1\,,
\end{CD}\ \ \ \ \\ \ \ \ \ \ \ 
\ee
$1\leq i \leq k-1$, 
imply that 
$$P_i^{-1}(J(G_i))=J(G_{i+1}), \ \ 1\leq i \leq k-1, \ \  P_k^{-1}(J(G_k))=J(G_{1}).$$ Thus, equality \eqref{xir} implies that 
$J(G_1)$ is a completely invariant set for $P_i,$ $1\leq i \leq k.$ In turn, this implies that 
$$J(P_1)=J(P_2)=\dots =J(P_k)=J(G_1)$$ (see \cite{be2}, Lemma 8, or \cite{p1}, Theorem 4). This proves the ``if'' part. 

Finally, to prove the ``only if'' part, we observe that 
if all elements of $S$ share a measure of maximal entropy, then by \eqref{zx} for every $i,$ $2\leq i \leq k,$ there exist $t_i,s_i\geq 1$ and $r_i\geq 0$ such that 
$$P_1^{\circ t_i}=P_1^{\circ r_i}\circ P_i^{\circ s_i}, \ \ \ \ 2\leq i \leq k.$$
Therefore, for $K=t_2\dots t_k$, we have:
 \be \l{ii} P_1^{\circ K}=(P_1^{\circ r_2}\circ P_2^{\circ s_2})^{\circ K/t_2}=
(P_1^{\circ r_3}\circ P_3^{\circ s_3})^{\circ K/t_3}=\dots = (P_1^{\circ r_3}\circ P_k^{\circ s_k})^{\circ K/t_k},\ee 
implying \eqref{int}. \qed

\end{document}